\documentclass[11pt]{article}
\usepackage{enumerate}
\usepackage{amssymb,a4wide,latexsym,makeidx,epsfig,fleqn}
\usepackage{amsthm}
\usepackage{amsmath}
\usepackage{enumerate}
\newtheorem{theorem}{Theorem}[section]

\newtheorem{lemma}[theorem]{Lemma}

\newtheorem{corollary}[theorem]{Corollary}

\begin{document}
\textwidth 150mm \textheight 225mm
\title{Sufficient conditions for Hamilton-connected graphs in terms of (signless Laplacian) spectral radius
\thanks{ Supported by
the National Natural Science Foundation of China (No. 11171273)}}
\author{{Qiannan Zhou$^1$, Ligong Wang$^{2,}$\footnote{Corresponding author.} and Yong Lu$^3$}\\
{\small Department of Applied Mathematics, School of Science, Northwestern
Polytechnical University,}\\ {\small  Xi'an, Shaanxi 710072,
People's Republic
of China.}\\
{\small $^1$ E-mail: qnzhoumath@163.com}\\
{\small $^2$ E-mail: lgwangmath@163.com}\\
{\small $^3$ E-mail: luyong.gougou@163.com}\\}
\date{}
\maketitle
\begin{center}
\begin{minipage}{120mm}
\vskip 0.3cm
\begin{center}
{\small {\bf Abstract}}
\end{center}
{\small In this paper, we present some spectral sufficient conditions for a graph to be Hamilton-connected in terms
of the spectral radius or signless Laplacian spectral radius of the graph. Our results improve some previous work.

\vskip 0.1in \noindent {\bf Key Words}: \ Hamilton-connected, sufficient condition, spectral radius, signless Laplacian. \vskip
0.1in \noindent {\bf AMS Subject Classification (2010)}: \ 05C50, 05C45. }
\end{minipage}
\end{center}

\section{Introduction }

In this paper, we only consider simple graphs without loops and multiple edges. For terminology and notation not
defined but used, we refer the reader to \cite{BondyM2008}. Let $G$ be a connected graph with vertex set
$V(G)=\{v_{1},v_{2},\ldots,v_{n}\}$ and edge set $E(G)$. Write by $m=e(G)=|E(G)|$ the number of edges of graph $G$.
For each $v_{i}\in V(G)$, denote by $N_{G}(v_{i})$ the set of vertices adjacent to $v_{i}$ in $G$ and $d_{i}=d_{G}(v_{i})=|N_{G}(v_{i})|$ the degree of $v_{i}$. Moreover,
$N[v]=N_{G}[v]=N_{G}(v)\cup \{v\}$. Denote by $\delta=\delta(G)$ the minimum degree of $G$ and $\Delta=\Delta(G)$ the maximum degree of $G$. For
convenience, we use ($0^{x_{0}}, 1^{x_{1}},\ldots, k^{x_{k}},\ldots,\Delta^{x_{\Delta}}$) to denote the degree sequence of $G$, where $x_{k}$ is the number of vertices of degree $k$ in $G$.
We use
$G+H$ and $G\vee H$ to denote the disjoint union and the join of $G$ and $H$ respectively. The union of $k$ disjoint
copies of the same graph $G$ is denoted by $kG$.

Let $G$ be a graph. The adjacency matrix and degree diagonal matrix of $G$ are denoted by $A(G)$ and $D(G)$, respectively. The largest eigenvalue of $A(G)$, denoted by $\rho(G)$, is called to be the spectral radius of $G$. The matrix
$Q(G)=D(G)+A(G)$ is the signless Laplacian matrix of $G$. The largest eigenvalue of $Q(G)$, denoted by $q(G)$, is
called to be the signless Laplacian spectral radius of $G$.

Hamilton path is a path containing all vertices of $G$, and Hamilton cycle is a cycle containing all vertices of $G$. A graph is called to be traceable if it contains a Hamilton path, and a graph is called to be Hamiltonian if it contains a Hamilton cycle. A graph is called to be Hamilton-connected if every two vertices of $G$ are connected by a Hamilton path.

The problem of determining whether a given graph is Hamiltonian, traceable, Hamilton-connected is NP-complete.
Recently, there are many reasonable sufficient or necessary conditions that were given for a
graph to be Hamiltonian, traceable or Hamilton-connected. Fiedler and Nikiforov \cite{Fiedler2010} firstly gave
sufficient conditions in terms of the spectral radius of a graph or its complement for the existence of Hamilton
cycles. This work motivated further research, one may refer to \cite{Benediktovich2015,Liu2015,Lu2012,NingGe2015,fy,YuYeCai2014,ZhouBo2010,ZhouQN22017}.
Recently, by imposing the minimum degree of a graph as a new parameter, Li and Ning \cite{LiNing2016LMA,LiNing2017LAA} extended some the results in \cite{Fiedler2010,Liu2015,NingGe2015}. Now, their
results were improved by Nikiforov \cite{Nikiforov2016}, Chen et al. \cite{ChenHou2017}, Ge et al. \cite{GeNing} and Li et al. \cite{LiLiuPeng}, in some sense.

The following sufficient condition involving the number of edge is due to Ore \cite{Ore1963}.

\noindent\begin{theorem}\label{th:5c2}(\cite{Ore1963}) Let $G$ be a graph on $n$ vertices and $m$ edges. If
$$m\geq \dbinom{n-1}{2}+3,$$
then $G$ is Hamilton-connected.
\end{theorem}

Observing that $\delta\geq 3$ is a trivial necessary condition for $G$ to be Hamilton-connected. Zhou and Wang
\cite{ZhouQN12017} refined the above edge number condition.

\noindent\begin{theorem}\label{th:5c3}(\cite{ZhouQN12017})
Let $G$ be a connected graph on $n\geq 6$ vertices and $m$ edges with minimum degree $\delta\geq 3$. If $$m\geq
\dbinom{n-2}{2}+6,$$ then $G$ is Hamilton-connected unless $G\in \mathbb{NP}_{1}=\{K_{3}\vee (K_{n-5}+2K_{1}),
K_{6}\vee 6K_{1}, K_{4}\vee (K_{2}+3K_{1}), 5K_{1}\vee K_{5}, K_{4}\vee (K_{1,4}+K_{1}), K_{4}\vee (K_{1,3}+K_{2}),
K_{3}\vee K_{2,5}, K_{4}\vee 4K_{1}, K_{3}\vee (K_{1}+K_{1,3}), K_{3}\vee (K_{1,2}+K_{2}), K_{2}\vee K_{2,4}\}.$
\end{theorem}

For $n\geq 5$ and $1\leq k\leq n/2$, we define:
$$S_{n}^{k}=K_{k}\vee (K_{n-(2k-1)}+(k-1)K_{1})~and~T_{n}^{k}=K_{2}\vee (K_{n-(k+1)}+K_{k-1}).$$
Moreover, for $t\geq 1$, let $\mathcal{S}_{n}^{k}(t)$ $(resp., \mathcal{T}_{n}^{k}(t))$ denote the set of all possible graphs obtained from $S_{n}^{k}$ ($resp., T_{n}^{k}$) by deleting exactly $t$ edges such that $\delta\geq 3$.
Obviously, $\mathcal{S}_{n}^{k}(0)=\{S_{n}^{k}\}$, $\mathcal{T}_{n}^{k}(0)=\{T_{n}^{k}\}$.

In this paper, we first make a further improvement for Theorem \ref{th:5c3}.

\noindent\begin{theorem}\label{th:5main1} Let $G$ be a connected graph on $n\geq 11$ vertices and $m$ edges with minimum degree $\delta\geq 3$. If
$$m\geq \dbinom{n-3}{2}+13,$$
then $G$ is Hamilton-connected unless $G\in (\bigcup_{i=0}^{n-10}\mathcal{S}_{n}^{3}(i))~\bigcup~
(\bigcup_{i=0}^{n-11}\mathcal{T}_{n}^{3}(i))$, or for $n=11$, $G=S_{11}^{5}$, or for $n=12$,
$G\in \bigcup_{i=0}^{2}\mathcal{S}_{12}^{6}(i)$, or for $n=13$, $G=S_{13}^{6}$, or for $n=14$,
$G\in \bigcup_{i=0}^{2}\mathcal{S}_{14}^{7}(i)~\bigcup~S_{14}^{3}(4)$,
 or for $n=16$,
$G\in \bigcup_{i=0}^{1}\mathcal{S}_{16}^{8}(i)~\bigcup~K_{7}\vee (K_{2}+K_{1,6})$.
\end{theorem}

By Theorem \ref{th:5main1}, we can get the following two corollaries immediately.

\noindent\begin{corollary}\label{co:5main2} Let $G$ be a connected graph on $n\geq 14$ vertices and $m$ edges with minimum degree $\delta\geq 3$. If
$$m\geq \dbinom{n-2}{2}+4,$$
then $G$ is Hamilton-connected unless $G\in (\bigcup_{i=0}^{2}\mathcal{S}_{n}^{3}(i))~\bigcup~
(\bigcup_{i=0}^{1}\mathcal{T}_{n}^{3}(i))$, or for $n=14$, $G=S_{14}^{7}$.

\end{corollary}

\noindent\begin{corollary}\label{co:5main3} Let $G$ be a connected graph on $n\geq 13$ vertices and $m$ edges with minimum degree $\delta\geq 3$. If
$$m\geq \dbinom{n-2}{2}+3,$$
then $G$ is Hamilton-connected unless $G\in (\bigcup_{i=0}^{3}\mathcal{S}_{n}^{3}(i))~\bigcup~
(\bigcup_{i=0}^{2}\mathcal{T}_{n}^{3}(i))$,
 or for $n=13$, $G=S_{13}^{6}$,
or for $n=14$, $G\in
\bigcup_{i=0}^{1}\mathcal{S}_{14}^{7}(i)$.
\end{corollary}

In \cite{fy}, Yu and Fan have established sufficient conditions for a graph to be Hamilton-connected in terms of the spectral radius and signless Laplacian spectral radius. Write $G=K_{n-1}+e+e'$ for $K_{n-1}$ together
with a vertex joining two vertices of $K_{n-1}$ by the edges $e$, $e'$, respectively.

\noindent\begin{theorem}\label{th:5cfy}(\cite{fy}) Let $G$ be a graph on $n$ vertices.
\begin{enumerate}[(i)]
  \item If $\rho(G)> -\frac{1}{2}+\sqrt{n^{2}-3n+\frac{17}{4}}$, then $G$ is Hamilton-connected unless $G=K_{n-1}+e+e'$.
  \item If $\rho(\overline{G})< \sqrt{\frac{(n-2)^{2}}{n}}$, then $G$ is Hamilton-connected.
  \item If $q(G)> 2n-4+\frac{2}{n-1}$, then $G$ is Hamilton-connected unless $G=K_{n-1}+e+e'$.
\end{enumerate}
\end{theorem}

Recently, Zhou and Wang \cite{ZhouQN12017} gave some spectral sufficient conditions on spectral radius and signless
Laplacian spectral radius for a graph to be Hamilton-connected, which extended the result of Yu and Fan \cite{fy} in
some sence.

\noindent\begin{theorem}\label{th:5czw}(\cite{ZhouQN12017}) Let $G$ be a connected graph on $n\geq 6$ vertices with minimum degree $\delta\geq 3$.
\begin{enumerate}[(i)]
  \item If $\rho(G)\geq \sqrt{n^{2}-6n+19}$, then $G$ is Hamilton-connected.
  \item If $q(G)\geq 2n-6+\frac{14}{n-1}$, then $G$ is Hamilton-connected unless $G=K_{4}\vee 4K_{1}$.
\end{enumerate}
\end{theorem}

In this paper, we continue to study new sufficient spectral conditions for a graph to be Hamilton-connected. We will
use Corollaries \ref{co:5main2} and \ref{co:5main3} to give the spectral sufficient conditions for a graph to be
Hamilton-connected.

\noindent\begin{theorem}\label{th:5main4} Let $G$ be a connected graph on $n\geq 14$ vertices with minimum degree $\delta\geq 3$. If
$$\rho(G)> n-3,$$
then $G$ is Hamilton-connected unless $G\in \{S_{n}^{3},T_{n}^{3}\}$.
\end{theorem}

It is easy to see that if we do 1 Kelmans operation on $T_{n}^{3}$, then we can obtain a proper
subgraph of $S_{n}^{3}$. Hence $\rho(S_{n}^{3})> \rho(T_{n}^{3})>\rho(K_{n-2})=n-3$. By Theorem \ref{th:5main4},
we have the following corollary.

\noindent\begin{corollary}\label{co:5main5} Let $G$ be a connected graph on $n\geq 14$ vertices with minimum degree $\delta\geq 3$. If
$$\rho(G)\geq \rho(S_{n}^{3}),$$
then $G$ is Hamilton-connected unless $G=S_{n}^{3}$.
\end{corollary}

\noindent\begin{theorem}\label{th:5main6} Let $G$ be a connected graph on $n\geq 13$ vertices with minimum degree $\delta\geq 3$. If
$$q(G)> 2n-6+\frac{6}{n-1},$$
then $G$ is Hamilton-connected unless $G=S_{n}^{3}$.
\end{theorem}

Obviously, we have the following corollary.

\noindent\begin{corollary}\label{co:5main7} Let $G$ be a connected graph on $n\geq 13$ vertices with minimum degree $\delta\geq 3$. If
$$q(G)\geq q(S_{n}^{3}),$$
then $G$ is Hamilton-connected unless $G=S_{n}^{3}$.
\end{corollary}

We can see that our results improve the previous work. Furthermore, let $\mathcal{G}_{n}$ be the class of
non-Hamilton-connected graphs of order $n$. In Corollaries \ref{co:5main5} and \ref{co:5main7}, we determine the maximum
spectral radius and the maximum signless
Laplacian spectral radius in $\mathcal{G}_{n}$. And the extremal graphs with maximum spectral radius and the maximum
signless Laplacian spectral radius are determined.

The rest of this paper is organized as follows. In Section~2, we present some useful techniques and lemmas. In
Section~3, we present the proofs of Theorems \ref{th:5main1},~\ref{th:5main4} and \ref{th:5main6}.

\section{Preliminaries}

In this section, we list some useful techniques and lemmas that will be used in later sections.

Firstly, let us recall the Kelmans transformation \cite{Kelmans1981}. Given a graph $G$ and two specified vertices
$u$, $v$ construct a new graph $G^{*}$ by replacing all edges $vx$ by $ux$ for $x\in N(v)\setminus N[u]$. Obviously,
the new graph $G^{*}$ has the same number of vertices and edges as $G$, and all vertices different from $u$ and $v$
retain their degrees. The vertices $u$ and $v$ are adjacent in $G^{*}$ if and only if they are adjacent in $G$.

\noindent\begin{lemma}\label{le:5c4}(\cite{Csikvari2009}) Let $G$ be a graph and $G^{*}$ be
a graph obtained from $G$ by some Kelmans's transformation. Then $\rho(G)\leq \rho(G^{*})$.
\end{lemma}

\noindent\begin{lemma}\label{le:5c5}(\cite{LiNing2016LMA}) Let $G$ be a graph and $G^{*}$ be
a graph obtained from $G$ by a Kelmans transformation. Then $q(G)\leq q(G^{*})$.
\end{lemma}

Suppose $M$ is a symmetric real matrix whose rows and columns are indexed by $X=\{1,\ldots,n\}$. Let
$\pi=\{X_{1},\ldots,X_{m}\}$ be a partition of $X$. Let $M$ be partitioned according to $\{X_{1},\ldots,X_{m}\}$,
i.e.,
\begin{displaymath}
A(\Gamma_{1})=\left(
  \begin{array}{cccccc}
        M_{11}&        \ldots &         M_{1m} \\
        \vdots&               &         \vdots \\
        M_{m1}&        \ldots &         M_{mm} \\

  \end{array}
\right),
\end{displaymath}
where $M_{ij}$ denotes the block of $M$ formed by rows in $X_{i}$ and the columns in $X_{j}$. Let $b_{ij}$ denote
the average row sum of $M_{ij}$, i.e., $b_{ij}=\frac{\textbf{1}^{T}M_{ij}\textbf{1}}{|X_{i}|}$, where $\textbf{1}$
is a column vector with all the elements 1. Then the matrix $M/\pi=(b_{ij})_{m\times m}$ is called the quotient
matrix of $M$. If the row sum of each block $M_{ij}$ is a constant, then the partition is called equitable.

\noindent\begin{lemma}\label{le:5c6}(\cite{GodsilRo2001}) Let $G$ be a graph. If $\pi$ is an equitable partition of
$V(G)$ corresponding to $A(G)$ $(Q(G))$, then $\rho(A(G)/\pi)=\rho(A(G))$ $(q(Q(G)/\pi)=q(Q(G)))$.
\end{lemma}

\noindent\begin{lemma}\label{le:5c7}(\cite{Brouwer2011,GodsilRo2001}) Let $G$ be a connected graph. If $H$ is a subgraph (proper subgraph) of $G$, then $\rho(H)\leq \rho(G)$ $(\rho(H)< \rho(G))$ and $q(H)\leq q(G)$ $(q(H)< q(G))$.
\end{lemma}

Hong et al. \cite{HongShu2001} proved the following spectral inequality for connected graphs. Nikiforov
\cite{Nikiforov2002} proved it for general graphs independently, and the case of equality was characterized in
\cite{ZhouCho2005}.

\noindent\begin{lemma}\label{le:5c8}(\cite{Nikiforov2002}) Let $G$ be a graph on $n$ vertices and $m$ edges with minimum degree $\delta$. Then $\rho(G)\leq
\frac{\delta-1}{2}+\sqrt{2m-n\delta+\frac{(\delta+1)^{2}}{4}}$.
\end{lemma}

The following result is also useful for us.

\noindent\begin{lemma}\label{le:5c9}(\cite{HongShu2001, Nikiforov2002}) For nonnegative integers $p$ and $q$ with
$2q\leq p(p-1)$ and $0\leq x\leq p-1$, the function $f(x)=\frac{x-1}{2}+\sqrt{2q-px+\frac{(1+x)^{2}}{4}}$
is decreasing with respect to $x$.
\end{lemma}

\noindent\begin{lemma}\label{le:5c10}(\cite{LFGU,fy}) Let $G$ be a connected graph on $n$ vertices and $m$ edges.
Then $q(G)\leq \frac{2m}{n-1}+n-2$.
\end{lemma}

\noindent\begin{lemma}\label{le:5c1}
(\cite{Berge,rao}) Let $G$ be a
graph on $n\geq 3$ vertices with degree sequence $(d_{1},d_{2},\ldots,d_{n})$, where $d_{1}\leq d_{2}\leq \cdots
\leq d_{n}$. If there is no integer $2\leq k\leq \displaystyle\frac{n}{2}$ such that $d_{k-1}\leq k$ and $d_{n-k}\leq
n-k$, then $G$ is Hamilton-connected.
\end{lemma}

\section{Proofs}



\noindent {\bf \emph{The proof of Theorem~\ref{th:5main1}}.} In this proof, we assume that a sequence $\vec{d}$ is called a
permissible graphic sequence
if there is a simple graph with degree sequence $\vec{d}$ satisfying the condition of Lemma \ref{le:5c1}. Suppose by contradiction that $G\notin (\bigcup_{i=0}^{n-10}\mathcal{S}_{n}^{3}(i))~\bigcup~
(\bigcup_{i=0}^{n-11}\mathcal{T}_{n}^{3}(i))$, and for $n=11$, $G\neq S_{11}^{5}$, and for $n=12$,
$G\notin \bigcup_{i=0}^{2}\mathcal{S}_{12}^{6}(i)$, and for $n=13$, $G\neq S_{13}^{6}$, and for $n=14$,
$G\notin \bigcup_{i=0}^{2}\mathcal{S}_{14}^{7}(i)~\bigcup~S_{14}^{3}(4)$, and for $n=16$,
$G\notin \bigcup_{i=0}^{1}\mathcal{S}_{16}^{8}(i)~\bigcup~K_{7}\vee (K_{2}+K_{1,6})$, but $G$ is
non-Hamilton-connected. Suppose that $G$ has the degree sequence $(d_{1},d_{2},\ldots,d_{n})$, where $d_{1}\leq
d_{2}\leq \cdots \leq d_{n}$. By Lemma \ref{le:5c1}, there exists an integer $k\leq n/2$ such that $d_{k-1}\leq k$ and $d_{n-k}\leq
n-k$. For convenience, we call this condition to be NHC-condition.  Thus, we have
\begin{align}\label{eq:1}
m&=\nonumber\frac{1}{2}\sum_{i=1}^{n}d_{i}\\
 &=\nonumber\frac{1}{2}(\sum_{i=1}^{k-1}d_{i}+\sum_{i=k}^{n-k}d_{i}+\sum_{i=n-k+1}^{n}d_{i})\\
 &\nonumber\leq \frac{1}{2}(k(k-1)+(n-k)(n-2k+1)+(n-1)k)\\
 &=\dbinom{n-3}{2}+13+\frac{f(k)}{2},
\end{align}
where $f(x):=3x^{2}-(2n+3)x+8n-38$. Since $e(G)\geq \dbinom{n-3}{2}+13$, combining with \eqref{eq:1}, we have
$f(k)>0$. Moreover, note that $3\leq \delta\leq d_{k-1}\leq k\leq n/2$, by a direct
computation, we obtain:
\begin{itemize}
  \item for $n\geq 11$, we have $f(3)=2n-20>0$, $f(4)=-2<0$.
\end{itemize}
Then we calculate $f(k)$ for $k\geq 5$. We have:

\begin{itemize}
  \item if $n=11$, then $k\leq 5$ and, $f(5)=0$;
  \item if $n=12$, then $k\leq 6$ and, $f(5)=-2<0$, $f(6)=4>0$;
  \item if $n=13$, then $k\leq 6$ and, $f(5)=-4<0$, $f(6)=0$;
  \item if $n=14$, then $k\leq 7$ and, $f(5)=-6<0$, $f(6)=-4<0$, $f(7)=4>0$;
  \item if $n=15$, then $k\leq 7$ and, $f(5)=-8<0$, $f(6)=-8<0$, $f(7)=-2<0$;
  \item if $n=16$, then $k\leq 8$ and, $f(5)=-10<0$, $f(6)=-12<0$, $f(7)=-8<0$, $f(8)=2>0$;
  \item if $n=17$, then $k\leq 8$ and, $f(5)=-12<0$, $f(6)=-16<0$, $f(7)=-14<0$, $f(8)=-6<0$;
  \item if $n\geq 18$, then for $5\leq k\leq n/2$, we have $f(k)<0$. To see this, we consider two roots of $f(x)=0$, which are
      $$r_{1}=\frac{2n+3-\sqrt{(2n+3)^{2}-12(8n-38)}}{6}, r_{2}=\frac{2n+3+\sqrt{(2n+3)^{2}-12(8n-38)}}{6}.$$
      By simple calculation, we have both $r_{1}<5$ and $r_{2}>n/2$ hold for $n\geq 18$, and then the desired result follows.
\end{itemize}

From the above computing results, we discuss the following cases.

\noindent \textbf{Case~1.} $k=3$ and $n\geq 11$.

In this case, we shall show when $G$ is not Hamilton-connected, $G\in (\bigcup_{i=0}^{n-10}\mathcal{S}_{n}^{3}(i))~\bigcup~
(\bigcup_{i=0}^{n-11}\\\mathcal{T}_{n}^{3}(i))$, which is a contradiction to our assumption.

By NHC-condition, we have
\begin{equation}\label{eq:2}
d_{1}=d_{2}=3, d_{3}\leq \cdots \leq d_{n-3}\leq n-3, d_{n-2}\leq d_{n-1}\leq d_{n}\leq n-1.
\end{equation}

Furthermore, note that $f(3)=2n-20$ when $n\geq 11$, by \eqref{eq:1}, we have
$$\dbinom{n-3}{2}+13\leq e(G)\leq \dbinom{n-3}{2}+13+(n-10).$$

If $e(G)=\dbinom{n-3}{2}+13+(n-10)$, then $\sum_{i=1}^{n}d_{i}=2e(G)=n^{2}-5n+18$. Hence it is direct that
the inequalities in \eqref{eq:2} must be equalities, and then the degree sequence of $G$ is
$(3,3,\underbrace{n-3,\ldots,n-3}_{n-5~times},n-1,n-1,n-1)$, which implies that
$G=S_{n}^{3}=\mathcal{S}_{n}^{3}(0)$, a contradiction.

If $e(G)=\dbinom{n-3}{2}+13+(n-10)-t$, where $1\leq t\leq n-10$, then we have
\begin{equation}\label{eq:3}
e(G)=e(S_{n}^{3})-t ~and ~ e(G)=e(T_{n}^{3})-(t-1).
\end{equation}

Moreover, note that any three 3-degree vertices are incident with at most 9 edges, and any $n-3$ vertices are
incident with $\dbinom{n-3}{2}$ edges, and $e(G)\geq \dbinom{n-3}{2}+13$, we conclude that $G$ has exactly two
3-degree vertices. Without loss of generality, we may suppose that $d_{G}(v_{1})=d_{G}(v_{2})=3$. Then we discuss
the following two subcases.

\noindent \textbf{Subcase~1.1.} $v_{1}$ is not adjacent to $v_{2}$.

If $N_{G}(v_{1})=N_{G}(v_{2})$, i.e., $v_{1}$ and $v_{2}$ have the same neighbour, then combining the definition
of $S_{n}^{k}$ and \eqref{eq:3}, $G$ is a subgraph of $S_{n}^{3}$ obtained by deleting $t$$(1\leq t\leq n-10)$
edges from its clique $K_{n-2}$. That is to say, $G\in \bigcup_{t=1}^{n-10}\mathcal{S}_{n}^{3}(t)$, which is a contradiction to our assumption.

Now we assume $N_{G}(v_{1})\neq N_{G}(v_{2})$. Let $H_{1}=G[V(G)\backslash \{v_{1}\}]$. Then, $|V(H_{1})|=n-1$,
$\delta(H_{1})\geq 3$ and $e(H_{1})=e(G)-3\geq \dbinom{n-3}{2}+10=\dbinom{(n-1)-2}{2}+10> \dbinom{(n-1)-2}{2}+6$.
Hence by Theorem \ref{th:5c3}, we have $H_{1}$ is Hamilton-connected. If every two vertices in $V(H_{1})$ can be connected by a
Hamilton path in $G$, then $G$ is also Hamilton-connected, a contradiction. Then there must exist two vertices
$w$ and $w'$ such that they are connected by a path passing through all vertices in $V(G)$ but not
$v_{1}$. Let $P$ be this path in a given direction (from $w$ to $w'$). Suppose the vertices
in $P$ are $w=y_{1},
y_{2}, \ldots, y_{n-1}=w'$ in sequence. Let $y_{i}$, $y_{j}$ and $y_{l}$ $(1\leq i <j <l\leq n-1)$ be three
vertices adjacent to $v_{1}$. Then it is obvious that $\{v_{1}, y_{i+1}, y_{j+1}\}$ is an independent set since $G$
is not Hamilton-connected. We claim that
\begin{equation}\label{eq:4}
d_{G}(y_{i+1})+d_{G}(y_{j+1})\leq n.
\end{equation}
To see this, consider the set $K=\{y_{r}|y_{r-1}\in N_{G}(y_{i+1})\cap V(H_{1}), r-1\leq i~or~r-1\geq j+2\}\bigcup
\{y_{s}|y_{s+1}\in N_{G}(y_{i+1})\cap V(H_{1}), i+1< s+1\leq j\}$. Note that $|K|\geq d_{G}(y_{i+1})-2$. This
follows
since the vertex $y_{i+1}$ is possibly adjacent to $y_{n-1}=w'$ and $y_{i+1}$ is both the successor of $y_{i}$
and the predecessor of $y_{i+2}$. Since $\{v_{1}, y_{i+1}, y_{j+1}\}$ is an independent set, we obtain
$d_{H_{1}}(y_{i+1})=d_{G}(y_{i+1})$, $d_{H_{1}}(y_{j+1})=d_{G}(y_{j+1})$, and $|K\cup N_{G}(y_{j+1})|\leq
|V(H_{1})|-|\{y_{j+1}\}|=n-2$. Thus, if $d_{G}(y_{i+1})+d_{G}(y_{j+1})\geq n+1$, then
\begin{align*}
|K\cap N_{G}(y_{j+1})|&=|K|+|N_{G}(y_{j+1})|-|K\cup N_{G}(y_{j+1})|\\
 &\geq d_{G}(y_{i+1})-2+d_{G}(y_{j+1})-(n-2)\\
 &\geq n+1-2-(n-2)=1,
\end{align*}
implying that $K$ and $N_{G}(y_{j+1})$ have a common vertex, say $y_{t}$. Obviously, $t\neq i+1, j+1$. If $t=i$,
then $y_{i+1}y_{j+1}\in E(G)$, a contradiction. If $t\leq i-1$, then $y_{i+1}y_{t-1}\in E(G)$. Then
$y_{1}\overrightarrow{P}y_{t-1}y_{i+1}\overrightarrow{P}y_{j}v_{1}y_{i}\overleftarrow{P}y_{t}y_{j+1}
\overrightarrow{P}y_{n-1}$
is a Hamilton path in $G$ connecting $y_{1}$ and $y_{n-1}$, a contradiction. If $t\geq j+2$, then $y_{i+1}y_{t-1}\in
E(G)$. Then $y_{1}\overrightarrow{P}y_{i}v_{1}y_{j}\overleftarrow{P}y_{i+1}y_{t-1}\overleftarrow{P}
y_{j+1}y_{t}\overrightarrow{P}y_{n-1}$ is a Hamilton path in $G$ connecting $y_{1}$ and $y_{n-1}$, a contradiction.
If $i+1< t\leq j$, then $y_{i+1}y_{t+1}\in E(G)$. Then $y_{1}\overrightarrow{P}y_{i}v_{1}y_{j}\overleftarrow{P}
y_{t+1}y_{i+1}\overrightarrow{P}y_{t}y_{j+1}\overrightarrow{P}y_{n-1}$ is a Hamilton path in $G$ connecting $y_{1}$
and $y_{n-1}$, a contradiction. Hence \eqref{eq:4} holds. Next, according the distribution of the neighbors of $v_{1}$,
we discuss the following subcases.

\noindent \textbf{Subcase~1.1.1.} We assume $y_{i+1}\neq y_{j-1}$ and $y_{j+1}\neq y_{l-1}$, then by the similar
method as above, we have

$$d_{G}(y_{j-1})+d_{G}(y_{l-1})\leq n.$$

Consequently, by considering the number of edges in $G$, we can get the following contradiction:

\begin{align*}
e(G)&\leq 3+(d_{G}(y_{i+1})+d_{G}(y_{j+1}))+(d_{G}(y_{j-1})+d_{G}(y_{l-1}))\\
&~~~+e(G[V(G)\setminus \{v_{1},y_{i+1},y_{j+1},y_{j-1},y_{l-1}\}])\\
 &\leq 3+n+n+\dbinom{n-5}{2}=\dbinom{n-3}{2}+12<\dbinom{n-3}{2}+13\\
 &\leq e(G),
\end{align*}
which is a contradiction and the first inequality follows from that there may be edges in vertex set $\{y_{i+1},y_{j+1},y_{j-1},y_{l-1}\}$.

\noindent \textbf{Subcase~1.1.2.} We assume $y_{i+1}= y_{j-1}$ and $y_{j+1}\neq y_{l-1}$, then $y_{j}=y_{i+2}$,
$y_{j+1}=y_{i+3}$.

If $y_{l}\neq y_{n-1}$, then by using the similar method as that of Subcase~1.1.1, we can obtain
\begin{align*}
e(G)&\leq 3+(d_{G}(y_{i+1})+d_{G}(y_{l-1}))+(d_{G}(y_{i+3})+d_{G}(y_{l+1}))\\
&~~~+e(G[V(G)\setminus \{v_{1},y_{i+1},y_{l-1},y_{i+3},y_{l+1}\}])\\
 &\leq 3+n+n+\dbinom{n-5}{2}=\dbinom{n-3}{2}+12<\dbinom{n-3}{2}+13\\
 &\leq e(G),
\end{align*}
which is a contradiction and the first inequality follows from that there may be edges in vertex set
$\{y_{i+1},y_{l-1},y_{i+3},y_{l+1}\}$.

Then we suppose $y_{l}=y_{n-1}$. If $y_{i}\neq y_{1}$, then by using the similar method as that of Subcase~1.1.1, we
can obtain
\begin{align*}
e(G)&\leq 3+(d_{G}(y_{i-1})+d_{G}(y_{l-1}))+(d_{G}(y_{i+1})+d_{G}(y_{i+3}))\\
&~~~+e(G[V(G)\setminus \{v_{1},y_{i-1},y_{l-1},y_{i+1},y_{i+3}\}])\\
 &\leq 3+n+n+\dbinom{n-5}{2}=\dbinom{n-3}{2}+12<\dbinom{n-3}{2}+13\\
 &\leq e(G),
\end{align*}
which is a contradiction and the first inequality follows from that there may be edges in vertex set $\{y_{i-1},y_{l-1},y_{i+1},y_{i+3}\}$.

If $y_{i}=y_{1}$, then $v_{1}$ is adjacent to $y_{1}$, $y_{3}$ and $y_{n-1}$. Let $W_{1}=V(G)\setminus \{v_{1},
y_{1}, y_{3}\}$. Then we show $\delta(G[W_{1}])\geq 3$.

If $d_{G[W_{1}]}(y_{2})=1$, then $y_{2}=v_{2}$. Since $N_{G}(v_{1})\neq N_{G}(v_{2})$ and $\{v_{1},y_{2},y_{4}\}$ is
an independent set, $W_{1}\setminus\{y_{4}, y_{n-1}\}$ has one vertex $w_{1}$ such that $y_{2}w_{1}\in E(G)$. Let
$W_{2}=W_{1}\setminus \{y_{2}\}$. Then
\begin{align*}
e(G[W_{2}])&\geq \dbinom{n-3}{2}+13-3-d_{G}(y_{2})-d_{G[W_{2}\cup \{y_{1},y_{3}\}]}(y_{1})-d_{G[W_{2}\cup \{y_{3}\}]}(y_{3})\\
&\geq \dbinom{n-3}{2}+13-3-3-(n-3)-(n-4)=\dbinom{(n-4)-1}{2}+5,
\end{align*}
which, together with the fact that $|W_{2}|=n-4$ and Theorem \ref{th:5c2}, we have $G[W_{2}]$ is
Hamilton-connected. Then there is a Hamilton path $w_{1}Py_{n-1}$ which connects $w_{1}$ and $y_{n-1}$ in
$G[W_{2}]$. Then $y_{1}v_{1}y_{3}y_{2}w_{1}Py_{n-1}$ is a Hamilton path connecting $y_{1}$ and $y_{n-1}$ in
$G$, a contradiction.

If $d_{W_{1}}(y_{2})=2$, then $d_{G}(y_{2})=4$. There always exists $w_{1}$ that we discussed above. Then
\begin{align*}
e(G[W_{2}])&\geq \dbinom{n-3}{2}+13-3-d_{G}(y_{2})-d_{G[W_{2}\cup \{y_{1},y_{3}\}]}(y_{1})-d_{G[W_{2}\cup \{y_{3}\}]}(y_{3})\\
&\geq \dbinom{n-3}{2}+13-3-4-(n-3)-(n-4)=\dbinom{(n-4)-1}{2}+4.
\end{align*}
Hence we can also get a contradiction by a similar method as above.

If $d_{W_{1}}(y_{4})=1$ or $2$, then $d_{G}(y_{4})=3$ or $4$. Let $W_{3}=W_{1}\setminus \{y_{4}\}$. By a similar
discussion as above, we can obtain $G[W_{3}]$ is Hamilton-connected and also get a contradiction.

If there exist a vertex $y_{k}\in W_{1}\setminus \{y_{2}, y_{4}\}$ which satisfies that $d_{G[W_{1}]}(y_{k})\leq 2$,
then $d_{G}(y_{k})\leq 4$. Since $d_{G}(y_{2})+d_{G}(y_{4})\leq n$,
\begin{align*}
e(G)&\leq d_{G}(v_{1})+d_{G}(y_{k})+(d_{G}(y_{2})+d_{G}(y_{4}))+e(G[V(G)\setminus\{v_{1},y_{2},y_{4},y_{k}\}])\\
&\leq 3+4+n+\dbinom{n-4}{2}=\dbinom{n-3}{2}+11<\dbinom{n-3}{2}+13,
\end{align*}
which is a contradiction.

Hence $\delta(G[W_{1}])\geq 3$. Note that $|W_{1}|=n-3$, and
\begin{align*}
e(G[W_{1}])&\geq \dbinom{n-3}{2}+13-3-d_{G[W_{1}\cup \{y_{1},y_{3}\}]}(y_{1})-d_{G[W_{1}\cup \{y_{3}\}]}(y_{3})\\
&\geq \dbinom{n-3}{2}+13-3-(n-2)-(n-3)=\dbinom{(n-3)-2}{2}+6.
\end{align*}
Then by Theorem \ref{th:5c3}, we have $G[W_{1}]$ is either Hamilton-connected, or $G[W_{1}]=K_{3}\vee (K_{(n-3)-5}+2K_{1})$.
If $G[W_{1}]$ is Hamilton-connected, then there is a Hamilton path connecting $y_{2}$ and $y_{n-1}$ in $G[W_{1}]$,
say $y_{2}Py_{n-1}$. Then $y_{1}v_{1}y_{3}y_{2}Py_{n-1}$ is a Hamilton path connecting $y_{1}$ and $y_{n-1}$ in $G$,
a contradiction. If $G[W_{1}]=K_{3}\vee (K_{(n-3)-5}+2K_{1})$, then $e(G[W_{1}])=\dbinom{(n-3)-2}{2}+6$,
$d_{G[W_{1}\cup \{y_{1},y_{3}\}]}(y_{1})=n-2$ and $d_{G[W_{1}\cup \{y_{3}\}]}(y_{3})=n-3$. Therefore, we have
$d_{G}(y_{1})=d_{G}(y_{3})=n-1$ and $G$ has only one 3-degree vertex $v_{1}$, which contradicts the fact that $G$
has exactly two 3-degree vertices.

Furthermore, the case of $y_{i+1}\neq y_{j-1}$ and $y_{j+1}= y_{l-1}$ can be proved in a similar method, thus we
omit it.

\noindent \textbf{Subcase~1.1.3.} We assume $y_{i+1}= y_{j-1}$ and $y_{j+1}= y_{l-1}$, then $y_{j}=y_{i+2}$,
$y_{l}=y_{i+4}$.

If $v_{1}$ is adjacent to neither $y_{1}$ nor $y_{n-1}$, then there must exist $y_{i-1}$ and $y_{i+5}$ since $n\geq
11$. Then by using a similar method as that of Subcase~1.1.1, we can obtain
\begin{align*}
e(G)&\leq 3+(d_{G}(y_{i-1})+d_{G}(y_{i+3}))+(d_{G}(y_{i+1})+d_{G}(y_{i+5}))\\
&~~~+e(G[V(G)\setminus \{v_{1},y_{i-1},y_{i+3},y_{i+1},y_{i+5}\}])\\
 &\leq 3+n+n+\dbinom{n-5}{2}=\dbinom{n-3}{2}+12<\dbinom{n-3}{2}+13\\
 &\leq e(G),
\end{align*}
which is a contradiction  and the first inequality follows from that there may be edges in vertex set $\{y_{i-1},y_{i+3},y_{i+1},y_{i+5}\}$.

If $v_{1}$ is adjacent to $y_{1}$, then $v_{1}$ is also adjacent to $y_{3}$ and $y_{5}$. One may easily get a
contradiction by a similar discussion as that of Subcase~1.1.2.

Similarly, if $v_{1}$ is adjacent to $y_{n}$, then $v_{1}$ is also adjacent to $y_{n-3}$ and $y_{n-5}$. One can also
get a contradiction by a similar discussion as that of Subcase~1.1.2.

\noindent \textbf{Subcase~1.2.} $v_{1}$ is adjacent to $v_{2}$.

Consider the graph $H_{2}:= G[V(G)\setminus \{v_{1},v_{2}\}]$. It is not difficult to see $|V(H_{2})|=n-2$
and $e(H_{2})=e(G)-5\geq \dbinom{(n-2)-1}{2}+8$. Then by Theorem \ref{th:5c2}, we get that $H_{2}$ is
Hamilton-connected. There must exist two vertices $w$ and $w'$ such that they are connected by a path
passing through all vertices in $V(H_{2})$ but not $v_{1}$ and $v_{2}$ at the same time. We denote this path by
$y_{1}P'y_{n-2}$, where $y_{1}=w$, $y_{n-2}=w'$, and give this path a direction (from $w$ to
$w'$). If $u$ is on this path, we use $u^{+}$ and $u^{-}$ to denote the successor and predecessor of $u$,
respectively.
 Since $d_{G}(v_{1})=d_{G}(v_{2})=3$, there must
be two vertices of $H_{2}$, say, $z_{1}$, $z_{2}$, (they are in order on this path)
which are adjacent to $v_{1}$. Also, there must be two vertices $z_{3}$ and $z_{4}$ (they are in order on this path) of $H_{2}$, which are adjacent
to $v_{2}$.

We now claim that $z_{1}=z_{3}$ and $z_{2}=z_{4}$, which, together with \eqref{eq:3}, would yield that $G$
is a subgraph of $T_{n}^{3}$ obtained by deleting $t-1$ edges from its clique $K_{n-2}$, that is, $G\in
\mathcal{T}_{n}^{3}(t-1)$, where $1\leq t\leq n-10$.

Suppose to the contrary that $z_{1}\neq z_{3}$ or $z_{2}\neq z_{4}$. We can easily see that $z_{i}\neq z_{1}^{-},z_{1}^{+},z_{2}^{-},z_{2}^{+}$ ($i=3,4$), $z_{j}\neq z_{3}^{-},z_{3}^{+},z_{4}^{-},z_{4}^{+}$ ($j=1,2$).
 And, if $z_{2}=z_{1}^{+}$,
then $z_{4}\neq z_{3}^{+}$ and vise versa. 

If $z_{1}\neq z_{3}$ and $z_{2}\neq z_{4}$, then by the same discussion on \eqref{eq:4}, we have
$d_{G}(z_{1}^{+})+d_{G}(z_{3}^{+})\leq n-1$. Then,
\begin{align*}
e(G)&=5+(d_{G}(z_{1}^{+})+d_{G}(z_{3}^{+}))+e(G[V(G)\setminus \{v_{1},v_{2},z_{1}^{+},z_{3}^{+}\}])\\
 &\leq 5+n-1+\dbinom{n-4}{2}=\dbinom{n-3}{2}+8<\dbinom{n-3}{2}+13\\
 &\leq e(G),
\end{align*}
which is a contradiction.

If $z_{1}= z_{3}$ and $z_{2}\neq z_{4}$, then by the same discussion on \eqref{eq:4}, we have
$d_{G}(z_{2}^{-})+d_{G}(z_{4}^{-})\leq n-1$. Then,
\begin{align*}
e(G)&=5+(d_{G}(z_{2}^{-})+d_{G}(z_{4}^{-}))+e(G[V(G)\setminus \{v_{1},v_{2},z_{2}^{-},z_{4}^{-}\}])\\
 &\leq 5+n-1+\dbinom{n-4}{2}=\dbinom{n-3}{2}+8<\dbinom{n-3}{2}+13\\
 &\leq e(G),
\end{align*}
which is a contradiction.
The case of $z_{1}\neq z_{3}$ and $z_{2}= z_{4}$ can be discussed in a similar way.

Summing up the above discussion, we have $G\in (\bigcup_{i=0}^{n-10}\mathcal{S}_{n}^{3}(i))~\bigcup~
(\bigcup_{i=0}^{n-11}\mathcal{T}_{n}^{3}(i))$, as desired.

\noindent \textbf{Case~2.} $n=11$ and $k=5$.

In this case, by NHC-condition, we have
\begin{equation}\label{eq:5}
d_{1}\leq d_{2}\leq d_{3}\leq d_{4}\leq5,~d_{5}\leq d_{6}\leq6,~d_{7}\leq\cdots\leq d_{11}\leq10.
\end{equation}
Moreover, since $f(5)=0$ when $n=11$, we obtain $e(G)=\dbinom{n-3}{2}+13=41$, and then $\sum_{i=1}^{11}d_{i}=82$.
Combining with \eqref{eq:5}, we have the degree sequence of $G$ is $(5^{4},6^{2},10^{5})$, which implies
$G=K_{5}\vee (K_{2}+4K_{1})=S_{11}^{5}$, a contradiction.

\noindent \textbf{Case~3.} $n=12$ and $k=6$.

Again, by NHC-condition, we have
\begin{equation}\label{eq:6}
d_{1}\leq d_{2}\leq d_{3}\leq d_{4}\leq d_{5}\leq 6, d_{6}\leq6,~d_{7}\leq\cdots\leq d_{12}\leq11.
\end{equation}
Note that $f(6)=4$ when $n=12$, and by \eqref{eq:1}, we have $49\leq e(G)\leq 51$.  If $e(G)=51$, then
$\sum_{i=1}^{12}d_{i}=102$, which, together with \eqref{eq:6}, yields that the degree sequence of $G$ is
$(6^{6},11^{6})$. From this one can check directly that $G=K_{6}\vee 6K_{1}=S_{12}^{6}= \mathcal{S}_{12}^{6}(0)$.

Now assume that
\begin{equation}\label{eq:7}
e(G)=51-t=e(S_{12}^{6})-t,~where~t\in \{1,2\}.
\end{equation}

Since $\sum_{i=1}^{12}d_{i}=102-2t\geq 98$, $G$ has at least one 6-degree vertex and has no 3-degree vertex. Let
$d_{G}(x_{0})=6$ and $H_{3}:=G[V(G)\setminus \{x_{0}\}]$. It is easy to see that $|V(H_{3})|=11$,
$\delta(H_{3})\geq 3$ and $e(H_{3})=e(G)-6\geq 49-6=43> \dbinom{11-2}{2}+6$, by Theorem \ref{th:5c3}, we have that $H_{3}$
is Hamilton-connected. Let $w P w'$ be a Hamilton path in $H_{3}$ from $w$ to $w'$. Since
$G$ is not Hamilton-connected, there is no Hamilton path connecting $w$ and $w'$ in $G$.
Suppose that $x_{1},x_{2},x_{3},x_{4},x_{5},x_{6}$ are the distinct vertices of $H_{3}$ which are adjacent to
$x_{0}$. Without loss of generality, we assume that $x_{6}= w'$. Then
$\{x_{0},x_{1}^{+},x_{2}^{+},x_{3}^{+},x_{4}^{+},x_{5}^{+}\}$ is an independent set, which, together with
\eqref{eq:7}, would yield that $G$ is a subgraph of $S_{12}^{6}$ obtained by deleting any $t$ edges, that is,
$G\in \mathcal{S}_{12}^{6}(t)$, where $t\in \{1,2\}$.

Summing up the above discussion, we eventually obtain $G\in \bigcup_{i=0}^{2}\mathcal{S}_{12}^{6}(i)$, a
contradiction.

\noindent \textbf{Case~4.} $n=13$ and $k=6$.

This case is completely analogous to Case ~2. We can obtain $e(G)=58$, $\sum_{i=1}^{13}d_{i}=116$, and the degree
sequence of $G$ is $(6^{5},7^{2},12^{6})$, which implies that $G=K_{6}\vee (K_{2}+5K_{1})=S_{13}^{6}$, a
contradiction.

\noindent \textbf{Case~5.} $n=14$ and $k=7$.

As Case~3, we have $d_{1}\leq \cdots\leq d_{6}\leq 7$, $d_{7}\leq 7$, $d_{8}\leq\cdots \leq d_{14}\leq 13$ and
$68\leq e(G)\leq 70$. Since $\sum_{i=1}^{14}d_{i}\geq 136$, $G$ has at least one 7-degree vertex and has no 3-degree
vertex. Let $d_{G}(x_{0})=7$ and $H_{4}=G[V(G)\setminus \{x_{0}\}]$. Obviously, $|V(H_{4})|=13$ and
$\delta(H_{4})\geq 3$.

If $e(G)=70$, then the degree sequence of $G$ must be $(7^{7}, 13^{7})$, which implies that $G=K_{7}\vee
7K_{1}=S_{14}^{7}= \mathcal{S}_{14}^{7}(0)$.

If $e(G)=70-1=69$, then $e(H_{4})=e(G)-7=62> \dbinom{13-2}{2}+6=61$, by Theorem \ref{th:5c3}, $H_{4}$ is
Hamilton-connected. Hence, by a similar argument as that in Case~3, we would get $G\in \mathcal{S}_{14}^{7}(1)$.

If $e(G)=70-2=68$, then $e(H_{4})=e(G)-7=61= \dbinom{13-2}{2}+6$, by Theorem \ref{th:5c3}, $H_{4}$ is either Hamilton-connected or
$H_{4}=K_{3}\vee (K_{8}+2K_{1})$. If $H_{4}$ is Hamilton-connected, by a similar argument as that in Case~3, we would get
$G\in \mathcal{S}_{14}^{7}(2)$, a contradiction. Hence, $H_{4}=K_{3}\vee (K_{8}+2K_{1})=S_{13}^{3}$. In this case, if $x_{0}$ is
not adjacent to the two 3-degree vertices of $H_{4}$, then it is evident that $G$ is a subgraph of $S_{14}^{3}$ with
$e(S_{14}^{3})-4$ edges, that is, $G\in \mathcal{S}_{14}^{3}(4)$; otherwise, one may check easily that $G$ is
Hamilton-connected, contradicting our assumption.

Summing up the above discussion, we eventually get $G\in \bigcup_{i=0}^{2}\mathcal{S}_{14}^{7}(i)\bigcup
\mathcal{S}_{14}^{3}(4)$.

\noindent \textbf{Case~6.} $n=16$ and $k=8$.

Similarly, we have $d_{1}\leq \cdots\leq d_{7}\leq 8$, $d_{8}\leq 8$, $d_{9}\leq \cdots\leq d_{16}\leq 15$,
$91\leq e(G)\leq 92$ and hence $182\leq \sum_{i=1}^{16}d_{i}\leq 184$. From the inequality
$\sum_{i=9}^{16}d_{i}=2m- \sum_{i=1}^{8}d_{i}\geq 182-64=118$, we get that $d_{11}=\cdots=d_{16}=15$ and
$d_{9}+d_{10}\geq 28$. Note that $\sum d_{i}$ is even and the total degree is between 182 and 184. If
$d_{9}=d_{10}=15$, then the permissible graphic sequence is $(8^{8},15^{8})$, which implies that
$G=K_{8}\vee 8K_{1}=S_{16}^{8}\in \mathcal{S}_{16}^{8}(0)$. If $d_{9}=14$ and $d_{10}=15$, then the permissible
graphic sequence is $(7^{1},8^{7},14^{1},15^{7})$, which implies $G=K_{7}\vee (K_{1}+K_{1,7})$. If $d_{9}=13$ and
$d_{10}=15$, then the permissible graphic sequence is $(8^{8},13^{1},15^{7})$, which implies that $G=K_{7}\vee
(K_{2}+K_{1,6})$. If $d_{9}=d_{10}=14$, then the permissible graphic sequence is $(8^{8},14^{2},15^{6})$. If $v_{9}$
is adjacent to $v_{10}$, then $G$ can be obtained as follows. Let $X=K_{8}$, $Y=8K_{1}$, $x_{1},x_{2}\in V (X)$ and
$y_{1},y_{2}\in V (Y)$. $G$ is a subgraph of $X\vee Y$ obtained by deleting $x_{1}y_{1}$, $x_{2}y_{2}$ and adding
a new edge $y_{1}y_{2}$. Note that in this case, $G$ is Hamilton-connected. If $v_{9}$ is not adjacent to $v_{10}$,
then $G=K_{6}\vee K_{2,8}$.

Since $\mathcal{S}_{16}^{8}(1)=\{K_{7}\vee (K_{1}+K_{1,7}), K_{6}\vee K_{2,8}\}$, summing up the above discussion,
we get $G\in \bigcup_{i=0}^{1}\mathcal{S}_{16}^{8}(i)\bigcup (G=K_{7}\vee (K_{2}+K_{1,6}))$.

The proof is complete.
\hfill$\blacksquare$






\noindent {\bf \emph{The proof of Theorem~\ref{th:5main4}}.} Suppose that $G$ is not Hamilton-connected. By Lemmas \ref{le:5c8} and \ref{le:5c9}, we have
$$\rho(G)\leq 1+\sqrt{2m-3n+4},$$
which, together with the condition of Theorem \ref{th:5main4}, yields
$$n-3< \rho(G)\leq 1+\sqrt{2m-3n+4}.$$
We obtain $2m> n^{2}-5n+12$. Furthermore, by parity, we have $m\geq \dbinom{n-2}{2}+4$. By Corollary
\ref{co:5main2}, we have $G\in (\bigcup_{i=0}^{2}\mathcal{S}_{n}^{3}(i))~\bigcup~
(\bigcup_{i=0}^{1}\mathcal{T}_{n}^{3}(i))$
or for $n=14$, $G=S_{14}^{7}$.

Note that $K_{n-2}$ is a proper subgraph of $S_{n}^{3}$ and $T_{n}^{3}$, by Lemma \ref{le:5c7}, we have
$\rho(S_{n}^{3})> \rho(K_{n-2})=n-3$ and $\rho(T_{n}^{3})> \rho(K_{n-2})=n-3$. So $S_{n}^{3}$ and $T_{n}^{3}$
enter the list of exceptions of the theorem.

For $G\in \mathcal{S}_{n}^{3}(1)$, that is, $G$ is obtained from the graph $S_{n}^{3}$ by removing one edge, which
can have only one of the following degree sequences:

\begin{enumerate}[(1)]
  \item $H_{1}$ has degree sequence $(3,3,\underbrace{n-3,\ldots,n-3}_{n-5~times},n-2,n-2,n-1)$, i.e., $H_{1}=K_{1,2}\vee (K_{n-5}+2K_{1})$;
  \item $H_{2}$ has degree sequence $(3,3,n-4,\underbrace{n-3,\ldots,n-3}_{n-6~times},n-2,n-1,n-1)$;
  \item $H_{3}$ has degree sequence $(3,3,n-4,n-4,\underbrace{n-3,\ldots,n-3}_{n-7~times},n-1,n-1,n-1)$, i.e., $H_{3}=K_{3}\vee ((2K_{1}\vee K_{n-7})+2K_{1})$.
\end{enumerate}

\begin{figure}[htbp]
  \centering
  \includegraphics[scale=0.6]{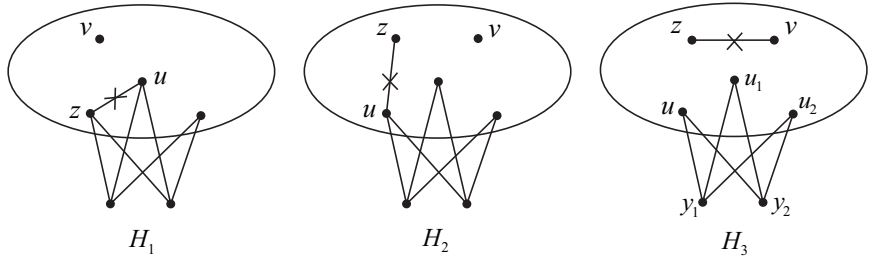}
\caption{Graphs, obtained from the graph $S_{n}^{3}$ by removing one edge.}
\end{figure}

The graphs which correspond to these degree sequences are depicted in Figure~1. Let $V_{1},V_{2}$ and $V_{3}$ be the
sets of vertices of $S_{n}^{3}$ with degree 3, $n-1$ and $n-3$. Therefore,
$H_{1}$ is the graph obtained from $S_{n}^{3}$ by deleting an edge $uz$ with $\{u, z\}\in V_{2}$. $H_{2}$ is the
graph obtained from $S_{n}^{3}$ by deleting an edge $uz$ with $u\in V_{2}$ and $z\in V_{3}$. $H_{3}$ is the graph
obtained from $S_{n}^{3}$ by deleting an edge $vz$ with $\{v,z\}\in V_{3}$.

Then we show $\rho(H_{i})< n-3$, $i=1,2,3$. Firstly, we claim that $\rho(H_{1})\leq \rho(H_{2})\leq \rho(H_{3})$.
Indeed, for graph $H_{i}$, let $u,z\in V(H_{i})$ be two vertices defined above and $v\in V_{3}\setminus\{u,z\}$,
where $i=1,2$. We have $H_{2}=H_{1}-vz+uz$ and $H_{3}=H_{2}-vz+uz$. Thus by Lemma \ref{le:5c4}, we can obtain the
conclusion. Hence it is sufficient to show only that $\rho(H_{3})< n-3$.

Let us consider the following partition of $V(H_{3})$ $\pi$: $X_{1}=V(K_{n-7})$, $X_{2}=\{z,v\}$, $X_{3}=\{u,
u_{1},u_{2}\}$, $X_{4}=\{y_{1},y_{2}\}$. It can easily be checked that this partition is equitable with the
adjacency matrix of its quotient $H_{3}/\pi$:

\begin{displaymath}
A(H_{3}/\pi)=\left(
  \begin{array}{cccccc}
        n-8&        2 &         3&     0 \\
         n-7&        0 &         3&     0\\
        n-7&        2 &         2&     2 \\
        0&        0 &         3&     0 \\

  \end{array}
\right).
\end{displaymath}

The characteristic polynomial $det(xI_{4}-A(H_{3}/\pi))$ of $A(H_{3}/\pi)$ is equal to:
\begin{equation}\label{eq:8}
f_{1}(x)=x^{4}-(n-6)x^{3}-(3n-7)x^{2}+4(n-10)x+12n-84.
\end{equation}

Therefore, by Lemma \ref{le:5c6}, the spectral radius of $H_{3}$ is the largest root of polynomial \eqref{eq:8}.

Next, we will show that there is no root of the  polynomial $f_{1}(x)$ in the interval $[n-3,+\infty)$. In fact,
when $n\geq 14$, it is obvious that the following inequalities are true:
\begin{align*}
f_{1}(n-3)&=2n^{2}-28n+18> 0;\\
f_{1}'(n-3)&=4x^{3}-3(n-6)x^{2}-2(3n-7)x+4(n-10)\mid_{x=n-3}\\
       &=n(n-3)^{2}-28> 0;\\
 f_{1}''(n-3)&=12x^{2}-6(n-6)x-2(3n-7)\mid_{x=n-3}=6n^{2}-24n+14>0;\\
 f_{1}'''(n-3)&=24x-6(n-6)\mid_{x=n-3}=18(n-2)>0;\\
 f_{1}^{(4)}(n-3)&=24>0.
\end{align*}

Therefore, by the Fourier-Budan theorem \cite{Prasolov2001}, all roots of $f_{1}(x)$ lie to the left of the number
$n-3$. In particular, $\rho(H_{3})< n-3$. Hence, non-Hamilton-connected graphs $G\in \mathcal{S}_{n}^{3}(1)$,
satisfy $\rho(G)< n-3$, a contradiction.

For $G\in \mathcal{S}_{n}^{3}(2)$, by Lemma \ref{le:5c7}, we also have $\rho(G)< n-3$, a contradiction.

For $G\in \mathcal{T}_{n}^{3}(1)$, that is, $G$ is obtained from the graph $T_{n}^{3}$ by removing one edge, which
can have only one of the following degree sequences:
\begin{enumerate}[(1)]
  \item $T_{1}$ has degree sequence $(3,3,\underbrace{n-3,\ldots,n-3}_{n-4~times},n-2,n-2)$, i.e., $T_{1}=2K_{1}\vee (K_{n-4}+K_{2})$;
  \item $T_{2}$ has degree sequence $(3,3,n-4,\underbrace{n-3,\ldots,n-3}_{n-5~times},n-2,n-1)$;
  \item $T_{3}$ has degree sequence $(3,3,n-4,n-4,\underbrace{n-3,\ldots,n-3}_{n-6~times},n-1,n-1)$, i.e., $T_{3}=K_{2}\vee ((2K_{1}\vee K_{n-6})+K_{2})$.
\end{enumerate}

\begin{figure}[htbp]
  \centering
  \includegraphics[scale=0.6]{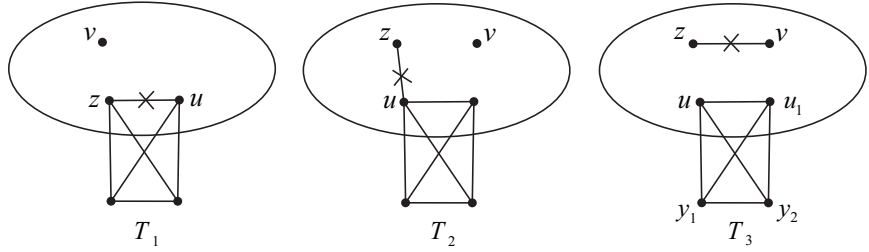}
\caption{Graphs, obtained from the graph $T_{n}^{3}$ by removing one edge.}
\end{figure}

The graphs which correspond to these degree sequences are depicted in Figure~2. Let $V_{1},V_{2}$ and $V_{3}$ be the
sets of vertices of $T_{n}^{3}$ with degree 3, $n-1$ and $n-3$. Therefore,
$T_{1}$ is the graph obtained from $T_{n}^{3}$ by deleting an edge $uz$ with $\{u, z\}\in V_{2}$. $T_{2}$ is the
graph obtained from $T_{n}^{3}$ by deleting an edge $uz$ with $u\in V_{2}$ and $z\in V_{3}$. $T_{3}$ is the graph
obtained from $T_{n}^{3}$ by deleting an edge $vz$ with $\{v,z\}\in V_{3}$.

Then we show $\rho(T_{i})< n-3$, $i=1,2,3$. Firstly, we claim that $\rho(T_{1})\leq \rho(T_{2})\leq \rho(T_{3})$.
Indeed, for graph $T_{i}$, let $u,z\in V(T_{i})$ be two vertices defined above and $v\in V_{3}\setminus\{u,z\}$,
where $i=1,2$. We have $T_{2}=T_{1}-vz+uz$ and $T_{3}=T_{2}-vz+uz$. Thus by Lemma \ref{le:5c4}, we can obtain the
conclusion. Hence it is sufficient to show only that $\rho(T_{3})< n-3$.

Let us consider the following partition of $V(T_{3})$ $\pi$: $X_{1}=V(K_{n-6})$, $X_{2}=\{z,v\}$, $X_{3}=\{u,
u_{1}\}$, $X_{4}=\{y_{1},y_{2}\}$. It can easily be checked that this partition is equitable with the
adjacency matrix of its quotient $T_{3}/\pi$:

\begin{displaymath}
A(T_{3}/\pi)=\left(
  \begin{array}{cccccc}
        n-7&        2 &         2&     0 \\
         n-6&        0 &         2&     0\\
        n-6&        2 &         1&     2 \\
        0&        0 &         2&     1 \\

  \end{array}
\right).
\end{displaymath}

The characteristic polynomial $det(xI_{4}-A(T_{3}/\pi))$ of $A(T_{3}/\pi)$ is equal to:
\begin{equation}\label{eq:9}
f_{2}(x)=x^{4}-(n-5)x^{3}-(2n-3)x^{2}+(5n-33)x+10n-56.
\end{equation}

Therefore, by Lemma \ref{le:5c6}, the spectral radius of $T_{3}$ is the largest root of polynomial \eqref{eq:9}.

Next, we will show that there is no root of the  polynomial $f_{2}(x)$ in the interval $[n-3,+\infty)$. In fact,
when $n\geq 14$, it is obvious that the following inequalities are true:
\begin{align*}
f_{2}(n-3)&=2n^{2}-20n+16> 0;\\
f_{2}'(n-3)&=4x^{3}-3(n-5)x^{2}-2(2n-3)x+5n-33\mid_{x=n-3}\\
       &=(n-3)^{2}(n-2)+n^{2}-7n-6> 0;\\
 f_{2}''(n-3)&=12x^{2}-6(n-5)x-2(2n-3)\mid_{x=n-3}=2(n-4)(3n-3)+2n>0;\\
 f_{2}'''(n-3)&=24x-6(n-5)\mid_{x=n-3}=6(3n-7)>0;\\
 f_{2}^{(4)}(n-3)&=24>0.
\end{align*}

Therefore, by the Fourier-Budan theorem \cite{Prasolov2001}, all roots of $f_{2}(x)$ lie to the left of the number
$n-3$. In particular, $\rho(T_{3})< n-3$. Hence, non-Hamilton-connected graphs $G\in \mathcal{T}_{n}^{3}(1)$,
satisfy $\rho(G)< n-3$, a contradiction.

For $n-14$, $G=S_{14}^{7}$, by direct calculation, we have $\rho(S_{14}^{7})=10.6158< 14-3=11$, a contradiction.

The proof is complete.
\hfill$\blacksquare$





\noindent {\bf \emph{The proof of Theorem~\ref{th:5main6}}.}  Combining Lemma \ref{le:5c10} and Theorem \ref{th:5main6}, we have
$$2n-6+\frac{6}{n-1}< q(G)\leq \frac{2m}{n-1}+n-2,$$
then $2m>n^{2}-5n+10$, which, by parity, is equivalent to $2m\geq n^{2}-5n+12=\dbinom{n-2}{2}+3$. Now, suppose that
$G$ is not Hamilton-connected, by Corollary \ref{co:5main3}, $G\in
(\bigcup_{i=0}^{3}\mathcal{S}_{n}^{3}(i))~\bigcup~
(\bigcup_{i=0}^{2}\\ \mathcal{T}_{n}^{3}(i))$,
or for $n=13$, $G=S_{13}^{6}$, or for $n=14$, $G\in
\bigcup_{i=0}^{1}S_{14}^{7}(i)$.

For $G=S_{n}^{3}$, it has been shown that $q(S_{n}^{3})$ is the largest zero of the function
$g_{1}(x)=x^{3}-(3n-5)x^{2}+(2n^{2}-n-24)x-6(n-3)(n-4)$ in \cite{ZhouQN12017}. Note that
$g_{1}(2n-6+\frac{6}{n-1})=-\frac{18(3n^{3}-18n^{2}+47n-44)}{(n-1)^{3}}< 0$ holds for $n\geq 13$, then we obtain
$$q(S_{n}^{3})>2n-6+\frac{6}{n-1},$$
so $S_{n}^{3}$ enter the list of exceptions of the theorem.

For $G\in \mathcal{S}_{n}^{3}(1)$, by the discussion in the proof of Theorem \ref{th:5main4} and Lemma \ref{le:5c5},
we have $q(H_{1})\leq q(H_{2})\leq q(H_{3})$. So it is sufficient to show only that $q(H_{3})<2n-6+\frac{6}{n-1}$.
Recall that the $(q,X)$-eigenequation in $G$ is
\begin{equation}\label{eq:10}
[q-d_{G}(v)]X_{v}=\sum_{u\in N_{G}(v)}X_{u},
\end{equation}
for each $v\in V(G)$, where $X$ is an eigenvalue of $Q(G)$ corresponding to the eigenvalue $q$ and $X_{u}$ is the
entry of $X$ corresponding to the vertex $u$. For $G=H_{3}$, let $X=(x_{1},x_{2},\ldots,x_{n})^{T}$ be the
eigenvector corresponding to $q(G)$. Then all vertices of degree $n-1$ have the same values given by
$X$, say $X_{1}$; all vertices of degree $n-3$ have the same values given by $X$, say $X_{2}$; all vertices of
degree $n-4$ have the same values given by $X$, say $X_{3}$.
Denote by $X_{4}$ the values of the vertices of degree 3 given by $X$. Assume
$\tilde{X}=(X_{1},X_{2},X_{3},X_{4})^{T}$, by \eqref{eq:10}, we have
\begin{align*}
(q(G)-(n-1))X_{1}&=2X_{1}+(n-7)X_{2}+2X_{3}+2X_{4};\\
(q(G)-(n-3))X_{2}&=3X_{1}+(n-8)X_{2}+2X_{3};\\
(q(G)-(n-4))X_{3}&=3X_{1}+(n-7)X_{2};\\
(q(G)-3)X_{4}&=3X_{1}.
\end{align*}

Transform the above equations into a matrix equation $(A-q(G)I)\tilde{X}=0$, we get
\begin{displaymath}
A=\left(
  \begin{array}{cccccc}
      n+1&        n-7&         2&        2&\\
        3&      2n-11&         2&        0&\\
        3&        n-7&       n-4&        0&\\
        3&          0&         0&        3&\\
  \end{array}
\right).
\end{displaymath}
Thus, $q(G)$ is the largest root of the following equation:
$$q^{4}-(4n-11)q^{3}+(5n^{2}-24n+10)q^{2}-(2n^{3}-7n^{2}-56n+220)q+6(n^{3}-13n^{2}+56n-80)=0.$$
Let $g_{2}(x)=x^{4}-(4n-11)x^{3}+(5n^{2}-24n+10)x^{2}-(2n^{3}-7n^{2}-56n+220)x+6(n^{3}-13n^{2}+56n-80)$, note that
$g_{2}(2n-6+\frac{6}{n-1})=\frac{2(4n^{6}-77n^{5}+445n^{4}-1471n^{3}+2939n^{2}-3856n+2664)}{(n-1)^{4}}> 0$ for
$n\geq 13$, which implies $q(H_{3})< 2n-6+\frac{6}{n-1}$. Hence $q(H_{1})\leq q(H_{2})\leq q(H_{3})<
2n-6+\frac{6}{n-1}$, a contradiction.

For $G\in \mathcal{S}_{n}^{3}(2)$, which is obtained from $S_{n}^{3}$ by deleting two edges, by Lemma \ref{le:5c7},
we also have $q(G)< 2n-6+\frac{6}{n-1}$, a contradiction.

For $G\in \mathcal{S}_{n}^{3}(3)$, which is obtained from $S_{n}^{3}$ by deleting three edges, by Lemma
\ref{le:5c7}, we also have $q(G)< 2n-6+\frac{6}{n-1}$, a contradiction.

For $G=T_{n}^{3}$, by a similar method as above, we get that the $q(G)$ is the largest zero of the function
$g_{3}(x)=x^{3}-(3n-4)x^{2}+2(n^{2}+n-14)x-8(n^{2}-6n+8)$. Note that
$g_{3}(2n-6+\frac{6}{n-1})=\frac{4(n^{4}-19n^{3}+102n^{2}-250n+220)}{(n-1)^{3}}> 0$ for $n\geq 13$, which implies
that $q(T_{n}^{3})< 2n-6+\frac{6}{n-1}$.

For $G\in \bigcup_{i=1}^{2}\mathcal{T}_{n}^{3}(i)$, which is a subgraph of $T_{n}^{3}$ by deleting $i\in \{1,2\}$
edges. By Lemma \ref{le:5c7}, we also have $q(G)< 2n-6+\frac{6}{n-1}$, a contradiction.

For $n=13$, $G=S_{13}^{6}$, by direct calculation, we have $q(S_{13}^{6})=20.1157< 2n-6+\frac{6}{n-1}$, a
contradiction.

For $n=14$, $G=S_{14}^{7}$, by direct calculation, we have $q(S_{14}^{7})=22.2195< 2n-6+\frac{6}{n-1}$,
which implies a contradiction. Hence, for $G\in \mathcal{S}_{14}^{7}(1)$, which is obtained from $S_{14}^{7}$ by
deleting one edge, by Lemma \ref{le:5c7}, we also have $q(G)< 2n-6+\frac{6}{n-1}$, a contradiction.



The proof is complete.
\hfill$\blacksquare$












\begin{thebibliography}{99}

\bibitem{Benediktovich2015} V. Benediktovich, Sufficient spectral condition for Hamiltonicity of a graph, Dokl. Nats. Akad. Nauk Belarusi, 59 (2015) 5--12.

\bibitem{Berge} C. Berge, Graphs and Hypergraphs. Translated from the French by Edward Minieka.
Second revised edition. North-Holland Mathematical Library, Vol. 6. North-Holland Publishing Co., Amsterdam-London; American Elsevier Publishing Co., Inc., New
York, 1976.

\bibitem{BondyM2008} J. A. Bondy, U. S. R. Murty, Graph Theory, Grad. Texts in Math, vol. 244, Springer, New York, 2008.

\bibitem{Brouwer2011} A. E. Brouwer, W. H. Haemaers, Spectra of Graphs, Springer-Verlag, 2011.

\bibitem{ChenHou2017} X. D. Chen, Y. P. Hou, J. G. Qian, Sufficient conditions for Hamiltonian graphs in terms of (signless Laplacian) spectral radius, Linear and Multilinear Algebra, http://dx.doi.org/10.1080/03081087.2017.1331994.

\bibitem{Csikvari2009} P. Csikvari, On a conjecture of V. Nikiforov, Discrete Math. 309 (2009) 4522--4526.

\bibitem{Erdos1962} P. Erd\H{o}s, Remarks on a paper of P\'{o}sa, Magyar Tud. Akad. Mat. Kut. Int. K\"{o}zl, 7 (1962) 227--229.

\bibitem{LFGU} L. H. Feng, G. H. Yu, On three conjectures involving the signless Laplacian spectral radius of graphs, Publ. Inst. Math. (Beograd) 85 (2009) 35--38.

\bibitem{Fiedler2010} M. Fiedler, V. Nikiforov, Spectral radius and Hamiltonicity of graphs, Linear Algebra Appl.
432 (2010) 2170--2173.

\bibitem{GeNing} J. Ge, B. Ning, Spectral radius and Hamiltonicity of graphs and balanced bipartite graphs with large minimum degree, Preprint available at arXiv:1606.08530v3.

\bibitem{GodsilRo2001} C. Godsil, C. Royle, Algebraic Graph Theory, Springer, 2001.

\bibitem{HongShu2001} Y. Hong , J. L. Shu, K. F. Fang, A sharp upper bound of the spectral radius of graphs, J. Comb. Theory B, 81 (2001) 177--183.

\bibitem{Kelmans1981} A. K. Kelmans, On graphs with randomly deleted edges, Acta Math. Acad. Sci. Hung. 37 (1981) 77--88.

\bibitem{LiNing2016LMA} B. L. Li, B. Ning, Spectral analogues of Erd\H{o}s' and Moon--Moser's theorems on Hamilton cycles, Linear and Multilinear Algebra, 64 (2016) 2252--2269.

\bibitem{LiNing2017LAA} B. L. Li, B. Ning, Spectral analogues of Moon--Moser's theorem on Hamilton paths in bipartite graphs, Linear Algebra Appl. 515 (2017) 180--195.

\bibitem{rao} R. Li, Laplacian spectral radius and some Hamiltonian properties of graphs,
    Appl. Math. E-Notes, 14 (2014) 216--220.

\bibitem{LiLiuPeng} Y. W. Li, Y. Liu, X. Peng, Signless Laplacian spectral radius and Hamiltonicity of graphs with large minimum degree, Linear Multilinear Algebra, http://dx.doi.org/10.1080/ 03081087.2017.1383346.

\bibitem{Liu2015} R. F. Liu, W. C. Shiu, J. Xue, Sufficient spectral conditions on Hamiltonian and traceable graphs,
Linear Algebra Appl. 467 (2015) 254--266.

\bibitem{Lu2012} M. Lu, H. Q. Liu, F. Tian, Spectral radius and Hamiltonian graphs, Linear Algebra Appl.
437 (2012) 1670--1674.

\bibitem{Nikiforov2002} V. Nikiforov, Some inequalities for the largest eigenvalue of a graph, Comb. Probab. Comput. 11 (2002) 177--183.

\bibitem{Nikiforov2011} V. Nikiforov, Some new results in extremal graph theory, Surveys in Combinatorics 2011, Cambridge University Press, 2011, 141--181.

\bibitem{Nikiforov2016} V. Nikiforov, Spectral radius and Hamiltonicity of graphs with large minimum degree, Czechoslovak Math. J. 66 (2016) 925--940.

\bibitem{NingGe2015} B. Ning, J. Ge, Spectral radius and Hamiltonian properties of graphs, Linear and Multilinear
Algebra, 63 (2015) 1520--1530.

\bibitem{Ore1963} O. Ore, Hamiltonian connected graphs, J Math Pures Appl., 42 (1963) 21--27.

\bibitem{Prasolov2001} V. V. Prasolov, Polynomials, MTSNMO, Moscow, 2001.

\bibitem{fy} G. D. Yu, Y. Z. Fan, Spectral conditions for a graph to be Hamilton-connected, Applied Mechanics and Materials, 336-338 (2013) 2329--2334.

\bibitem{YuYeCai2014} G. D. Yu, M. L. Ye, G. X. Cai, J. D. Cao, Signless Laplacian spectral conditions for Hamiltonicity of graphs, J. Appl. Math. 2014. Art. ID 282053, 6 PP.

\bibitem{ZhouBo2010} B. Zhou, Signless Laplacian spectral radius and Hamiltonicity, Linear Algebra Appl.
432 (2010) 566--570.

\bibitem{ZhouCho2005} B. Zhou, H. H. Cho, Remarks on spectral radius and Laplacian eigenvalues of a graph, Czechoslovak Math. J. 55 (2005) 781--790.

\bibitem{ZhouQN22017} Q. N. Zhou, L. G. Wang, Distance signless Laplacian spectral radius and Hamiltonian
properties of graphs, Linear Multilinear Algebra, http://dx.doi.org/10.1080/
03081087.2016.1273314.

\bibitem{ZhouQN12017} Q. N. Zhou, L. G. Wang, Some sufficient spectral conditions on Hamilton-connected and
traceable graphs, Linear Multilinear Algebra, 65 (2017) 224--234.






















































































\end{thebibliography}
\end{document}